\documentclass[10pt]{amsart}
\usepackage{amsmath,amssymb,amsthm}

\usepackage[pagebackref]{hyperref}
\usepackage{pinlabel}
\usepackage{overpic}
\usepackage{svg}
\usepackage{sseq}
\hypersetup{
    colorlinks = true,
    citecolor = [rgb]{0, 0.5, 0.00},
}

\renewcommand*{\backref}[1]{}
\renewcommand*{\backrefalt}[4]{%
    \ifcase #1 (Not cited.)%
    \or        (Cited on page~#2.)%
    \else      (Cited on pages~#2.)%
    \fi}

\usepackage{mathtools}
\usepackage{amssymb}

\usepackage[margin=1.175in]{geometry}
\usepackage{todonotes}
\usepackage{multirow}
\usepackage{longtable}
\usepackage{enumerate}
\usepackage{mathabx}
\usepackage{amsmath}
\usepackage{amsfonts}
\usepackage{amssymb}
\usepackage{float}
\usepackage{amsthm}
\usepackage{comment}
\usepackage{amscd}
\usepackage{amsxtra}
\usepackage{epsfig}
\usepackage{slashed}

\usepackage{listings}
\usepackage{color}

\definecolor{dkgreen}{rgb}{0,0.6,0}
\definecolor{gray}{rgb}{0.5,0.5,0.5}
\definecolor{mauve}{rgb}{0.58,0,0.82}

\usepackage{color}
\usepackage[all]{xy}
\usepackage{verbatim}

\usepackage{eucal}
\usepackage{mathrsfs}

\usepackage{graphicx}
\usepackage{tikz-cd}

\usepackage{url}
\usepackage{verbatim}

\usepackage{xy}
\xyoption{all}

\usepackage{amsthm}
\usepackage{tikz}
\usetikzlibrary{decorations.pathreplacing}

\usepackage{caption} 
\makeatletter
\def\@tocline#1#2#3#4#5#6#7{\relax
  \ifnum #1>\c@tocdepth 
  \else
    \par \addpenalty\@secpenalty\addvspace{#2}%
    \begingroup \hyphenpenalty\@M
    \@ifempty{#4}{%
      \@tempdima\csname r@tocindent\number#1\endcsname\relax
    }{%
      \@tempdima#4\relax
    }%
    \parindent\z@ \leftskip#3\relax \advance\leftskip\@tempdima\relax
    \rightskip\@pnumwidth plus4em \parfillskip-\@pnumwidth
    #5\leavevmode\hskip-\@tempdima
      \ifcase #1
       \or\or \hskip 1em \or \hskip 2em \else \hskip 3em \fi%
      #6\nobreak\relax
    \hfill\hbox to\@pnumwidth{\@tocpagenum{#7}}\par
    \nobreak
    \endgroup
  \fi}
\makeatother

\newtheorem{lemma}{Lemma}[section]

\newtheorem{corollary}[lemma]{Corollary}

\newtheorem{theorem}[lemma]{Theorem}

\newtheorem{prop}[lemma]{Proposition}

\theoremstyle{definition}

\newtheorem{definition}[lemma]{Definition}

\newtheorem{remark}[lemma]{Remark}

\numberwithin{equation}{section} 
\newcommand{\Z}{\mathbb{Z}}

\newcommand{\R}{\mathbb{R}}

\renewcommand{\ss}{\mathfrak{s}}

\newcommand{\Res}{\mathrm{Res}}

\usepackage{tikz-cd}

\newcommand{\HH}{\mathbb{H}}

\title{Vanishing Theorem of Pin(2)-equivariant Bauer-Furuta Invariants as Stable Sections}
\author{Haochen Qiu}
\begin{document}
\maketitle
\begin{abstract}
    We develop a viewpoint of $\operatorname{Pin}(2)$-equivariant Bauer-Furuta invariant: It can be converted to a section over a finite dimensional $\operatorname{Pin}(2)$-orbit, and it takes values in ordinary stable homotopy class. This viewpoint reveals a topological limitation of the invariant: If the fiber's connectivity exceeds the dimension of the base space, then all sections are homotopic to the trivial section. We use it to recover Lin-Mukherjee vanishing theorem. Furthermore, the family Seiberg-Witten invariant can be obtained from the pullback of such sections to $S^1$-orbit. We use this fact to prove a vanishing theorem of the $\operatorname{Pin}(2)$-equivariant family Bauer-Furuta invariant after two stabilizations.
\end{abstract}

\tableofcontents

\section{Introduction}

The Bauer-Furuta invariant \cite{BF02} as a smooth invariant, can be used to detect diffeomorphism of $4$-manifold which is not smoothly isotopic to the identity. By a finite-dimensional approximation, it is an equivariant stable mapping class for a spin manifold $X$:
\begin{equation*}\label{equ:group-of-BF}
\text{BF}^{\text{Pin}(2)} (X, \ss) \in \{ TF_0,  S^{n\mathbb{H}+b^+_2(X)\tilde{\R}} \}^{\text{Pin}(2)}.
\end{equation*}
where $\ss$ is a $\text{Spin}^c$-structure of $X$, and $TF_0$ is the Thom space of a rank $m$ quaternion bundle over $Pic^{\ss}(X) = T^{b_1(X)}$, such that 
\[
m-n = -\frac{\sigma(X)}{16}
\]
where $\sigma(X)$ is the signature of $X$. One can also forget the $\text{Pin}(2)$-action and define the nonequivariant Bauer-Furuta invariant by 
\[
\text{BF}^{\{e\}} (X, \ss):= \text{Res}^{\text{Pin}(2)}_{\{e\}} \text{BF}^{\text{Pin}(2)} (X, \ss) \in \{ TF_0,  S^{4n+b^+_2(X)}\}.
\]

For a spin family of simply-connect $4$-manifolds over $S^1$, one may define the family Bauer-Furuta invariant. We use Lin's convention~\cite[Definition 2.21]{Lin23}: 
\begin{equation}\label{equ:group-of-BF}
\text{FBF}^{\text{Pin}(2)} (X, \ss) \in \{ S^{\mathbb R} \wedge S^{m\mathbb H},  S^{n\mathbb{H}+b^+_2(X)\tilde{\R}} \}^{\text{Pin}(2)},
\end{equation}
where $S^{\mathbb R}$ comes from the parameter circle. Kronheimer-Mrowka~\cite{KM20} used the non-equivariant family Bauer-Furuta invariant to show that the Dehn twist along the separating $3$-sphere in the neck of $K3\mathbin{\#} K3$ is not smoothly isotopic to the identity. Lin~\cite{Lin23} used the $\operatorname{Pin}(2)$-equivariant family Bauer-Furuta invariant to show that, the above Dehn twist on $K3\mathbin{\#} K3$ is still not smoothly isotopic to the identity after a stabilization by $S^2 \times S^2$. The author~\cite{qiu25T4} used a refinement of $S^1$-equivariant family Bauer-Furuta invariant to show that, the Dehn twist along the separating $3$-sphere in the neck of any connect sum of odd homology $4$-tori is not smoothly isotopic to the identity.

However, there are some topological limitations of Bauer-Furuta invariant. Lin and Mukherjee~\cite{LM21} proved that, when $2\chi(X) + 3\sigma(X) \ge 4$, the equivariant Bauer-Furuta invariant vanishes for a large class of diffeomorphisms:
\begin{theorem}[\cite{LM21} Theorem 1.9]\label{thm:LM-vanishing}
    Let X be an oriented spin 4-manifold with $b_1(X) = 0$. We fix a homology orientation and
a spin structure on X. Let $f : X \to X$ be a diffeomorphism which preserves the spin structure and acts
trivially on $H^*(X; \mathbb R)$. Then for any lift $\widetilde{f}$ of $f$ to an automorphism on the spin structure, we have the
following results:
\begin{enumerate}[(1)]
 \item Suppose $2\chi(X) + 3\sigma(X) = 4$. Then the invariants
$\operatorname{FBF}^{S^1}(X, f, \widetilde{f} )$ vanishes.
 \item Suppose $2\chi(X) + 3\sigma(X) > 4$. Then the invariants
$\operatorname{FBF}^{\operatorname{Pin}(2)}(X, f, \widetilde{f} )$ vanishes.
\end{enumerate}
\end{theorem}
They first prove that the family BF invariant of such diffeomorphism must be so-called \textit{BF-type}, which means that the restriction of the Bauer-Furuta map to $S^1$-fixed point (i.e., remove all quaternion in both domain and target) is trivial. Then they use pure topological method including tom Dieck splitting theorem in equivariant stable homotopy theory, to show such restriction is injective. Therefore, the family BF invariant is not able to detect such diffeomorphism.

In this paper we investigate this topological limitation from another viewpoint. First we define
\begin{definition}
    For a compact space $Q$ and two finite dimensional vector bundle $V$ and $W$ over $Q$, we can define the section of the fiberwise stable mapping class space
    \[
    \{S^{V_q}, S^{W_q}\}=\operatorname*{colim}_{a\to\infty}
\operatorname{Map}_*(S^{\mathbb R^a\oplus V_q},S^{\mathbb R^a\oplus W_q}).
    \]
    We call such section a \textit{stable section}. It can be realized by a bundle map 
    \[
    V \oplus \underline{\mathbb R}^{a} \to W \oplus \underline{\mathbb R}^{a} 
    \]
    with $a$ sufficiently large, which only depends on $\operatorname{rank}(V)$ and $\operatorname{rank}(W)$. 
\end{definition}
The motivation is: For $S^1$-equivariant invariant, we can quotient out the domain by $S^1$-action, and convert it to a non-equivariant map from the $S^1$-orbit; However for the $\operatorname{Pin}(2)$-equivariant invariant, the group action acts on the target nontrivially, so it can only be converted to a fiberwise stable map over the $\operatorname{Pin}(2)$-orbit. In Proposition \ref{prop:G-equiv-to-stable-section}, we prove that the $\operatorname{Pin}(2)$-equivariant Bauer-Furuta map (\ref{equ:group-of-BF}) can be converted into a stable section (whose values are non-equivariant stable maps $\{S^2, S^{b^+}\}$) over the $\operatorname{Pin}(2)$-orbit of the domain. Then the topological limitation basically comes from the dimension of this base space: If the connectivity of the fiber is higher than this dimension, then there is no nontrivial section.

The connectivity of the fiber $\{S^{V_q}, S^{W_q}\}$ is $\operatorname{dim}W_q- \operatorname{dim}V_q$, and in the case of Bauer-Furuta map, it is given by $\approx b^+$ (adjusted by the dimension of the parameter space). The dimension of the base space is given by $-\frac{\sigma}{4}$. The difference between them is proportional to $4+5b^+-b^-$, which is exactly $2\chi(X) + 3\sigma(X)$ in Theorem \ref{thm:LM-vanishing}. Therefore, from this viewpoint, we can recover Theorem \ref{thm:LM-vanishing}.

Furthermore, the stabilization by $S^2\times S^2$ increases the connectivity of the fiber. Hence from the stable section viewpoint we can release the requirement on $2\chi(X) + 3\sigma(X)$ while proving the vanishing theorem of diffeomorphisms after stabilizations by $S^2\times S^2$:
\begin{remark}
$\Sigma^{n}P_+$ means $\Sigma^{n}(P_+)$. 
\end{remark}

\begin{theorem}\label{thm:simply-connected-two-euler-annihilation}
Let $G=\mathrm{Pin}(2)=S^1\cup jS^1\subset\mathbb H^\times$. Let $\tilde{\R}$ be its real sign representation, trivial on $S^1$, with $j$ acting as $-1$; let $\mathbb H$ have the left multiplication action; and let $\mathbb R$ be the trivial representation. Write $\{X,Y\}^{G}$ for genuine stable based equivariant maps in a complete universe, and write $e:S^0\to S^{\tilde{\R}}$ for the Euler inclusion. For an integer $n\ge2$, let $i_n:S^{\mathbb R}\to S^{\mathbb R+n\mathbb H}$ be inclusion at zero in the quaternionic coordinates. Define
\[
t_n:\{S^{\mathbb R+n\mathbb H},S^{(4n-2){\tilde{\R}}}\}^{S^1}\longrightarrow H^{4n-2}( \Sigma^2\mathbb{CP}^{2n-1}_+ ;\mathbb Z)\cong\mathbb  Z
\]
by the isomorphism with $\{\Sigma^2\mathbb{CP}^{2n-1}_+,S^{4n-2}\}$ induced by the cofiber sequence
\[
S^{\mathbb R}\longrightarrow S^{\mathbb R+n\mathbb H}\longrightarrow\Sigma^2S(n\mathbb H)_+,
\]
followed by the induced map on $H^{4n-2}(-; \mathbb Z)$ and pullback of $1\in H^{4n-2}(S^{4n-2}; \mathbb Z)$.

If $\alpha\in\{S^{\mathbb R+n\mathbb H},S^{(4n-2){\tilde{\R}}}\}^{G}$ satisfies
\[
\alpha\circ i_n=0,
\qquad
t_n(\operatorname{Res}_{S^1}^{G}\alpha)=0,
\]
then $\alpha e^2=0$.
\end{theorem}

This vanishing theorem shows that the $\operatorname{Pin}(2)$-equivariant FBF is not able to detect diffeomorphisms after two stabilizaitons by $S^2\times S^2$, for $4$-manifolds with $2\chi(X) + 3\sigma(X) \ge -4$ and homology-preserving diffeomorphisms. 
\begin{corollary}\label{coro:main-theorem-2stabilization-vanish}
    Let X be an oriented spin 4-manifold with $b_1(X) = 0$. We fix a homology orientation and
a spin structure on X. Let $f : X \to X$ be a diffeomorphism which preserves the spin structure and acts
trivially on $H^*(X; \mathbb R)$. Suppose $2\chi(X) + 3\sigma(X) \ge -4$. For any lift $\widetilde{f}$ of $f\# \operatorname{id}_{2(S^2\times S^2)}$ to an automorphism on the spin structure, if the family Seiberg-Witten invariant vanishes, then $\operatorname{FBF}^{\operatorname{Pin}(2)}(X, f\# \operatorname{id}_{2(S^2\times S^2)}, \widetilde{f} )$ vanishes.
\end{corollary}


In particular, in Lin's example~\cite{Lin23}, the $\operatorname{Pin}(2)$-equivariant family Bauer-Furuta invariant is zero for the Dehn twist on $K3\mathbin{\#} K3$ after two stabilizations by $S^2 \times S^2$.
\begin{corollary}\label{coro:lin-K3}
Let $\delta$ be the separating-neck Dehn twist on $K3\#K3$, and let $\alpha_0,\alpha_1\in\{S^{\mathbb R+2\mathbb H},S^{6{\tilde{\R}}}\}^{G}$ be the $\mathrm{Pin}(2)$-equivariant family Bauer-Furuta classes of its mapping torus with its two family spin structures defined in \cite{Lin23} Section 3.2. Then, for $s=0,1$,
\[
\alpha_s e\ne0,
\qquad
\alpha_s e^2=0.
\]
Consequently both family spin structures of the mapping torus after two identity $S^2\times S^2$ stabilizations have zero $\mathrm{Pin}(2)$-equivariant Bauer-Furuta invariant.
\end{corollary}

Of course, this statement can be deduced from Theorem \ref{thm:LM-vanishing}: $K3\#2(S^2\times S^2)$ has $2\chi(X) + 3\sigma(X)=8$, then use gluing formula of the family BF invariants to glue another $K3$.  But Corollary \ref{coro:main-theorem-2stabilization-vanish} also works for the diffeomorphisms supported on the the whole $4$-manifold, where the gluing formula does not apply.


\subsection*{Acknowledgements}
The author wants to thank Yi Sha for valuable conversations, especially for pointing out that the two-stabilization-vanishing theorem of the connected sum of $K3$ can be deduced from \cite{LM21} and the gluing formula.
\subsection*{AI Disclosure} The huge commutative diagram and the spectral sequence in Prop \ref{prop:obtain-SW-from-primary-obstruction}, and the picture in Lemma \ref{lem:alpha-i-trival}, are first hand-drawn and then converted to Latex code by ChatGPT 6 Astra. The originality of the work was checked by Codex-based knowledge navigator Qprint~\cite{github_hcqiu_qprint}.

\section{Convert family Bauer-Furuta and family Seiberg-Witten to stable sections}

\subsection{Convert the $G$-equivariant map to a stable section} In this subsection we convert the $G$-equivariant map to a stable section and consider the relevant primary obstruction. We adopt the settings in Theorem \ref{thm:simply-connected-two-euler-annihilation}. Set
\[
P=S(n\mathbb H),\qquad Q=S(n\mathbb H)/G,\qquad\lambda=S(n\mathbb H)\times_G{\tilde{\R}},
\]
where $\lambda$ is the associated real line bundle on $Q$. The action on $S(n\mathbb H)$ is free: if a unit quaternion fixes a nonzero component of a vector, that quaternion is $1$. The cofiber sequence
\[
S^{\mathbb R}\xrightarrow{i_n} S^{\mathbb R+n\mathbb H}\xrightarrow{c_n}
\Sigma^2S(n\mathbb H)_+
\]
and exactness give a class
\[
\beta\in\{\Sigma^2S(n\mathbb H)_+,S^{(4n-2){\tilde{\R}}}\}^{G}
\]
such that
\[
\alpha =\beta c_n.
\]
We do not assume that the lift is unique.

Here is the identification we will prove:
\[
\text{$G$-equivariant map from $S(n\mathbb H)$}
\quad\longleftrightarrow\quad
\text{Mapping class family over $S(n\mathbb H)/G$}.
\]
Before we start, we note that
\begin{lemma}\label{lem:based-map-adjunction} Suppose $P$ is a space and $V, W$ are two $\mathbb R$-vector spaces. Regard the infinity point of the one point compactification of $\mathbb R$-vector space to be the base point. Then
\[
\operatorname{Map}(P,\operatorname{Map}_*(S^V,S^W)) \cong \operatorname{Map}_*(\Sigma^V(P_+), S^W).
\]
\end{lemma}
This is because to obtain $\Sigma^V(P_+)=S^V \wedge (P_+)$ we take the product and then collapse a copy of $P$ at the infinity point of $S^V$. Hence for any point of $P$, a map in $\operatorname{Map}_*(\Sigma^V(P_+), S^W)$ gives a based map $S^V \to S^W$. The inverse is also straightforward.
\begin{prop}\label{prop:G-equiv-to-stable-section}
For a compact free $G$-space $P$, a class in $\{\Sigma^2P_+,S^{k{\tilde{\R}}}\}^{G}$ is a homotopy class of sections, over $Q=P/G$, of the fiberwise stable mapping class space
\[
\mathcal E_k(q)=\{S^2,S^{k\lambda_q}\} := \operatorname*{colim}_{a\to\infty}
\operatorname{Map}_*(S^{a+2},S^{\mathbb R^a\oplus k\lambda_q}),
\]
where $\lambda=P\times_G{\tilde{\R}}$ and $\lambda_q \cong {\tilde{\R}}$ is the fiber of $\lambda$ over $q\in Q$. The zero section is the constant map to infinity. 
\end{prop}
To see the description, first represent the stable map with only a sufficiently large trivial suspension on the $G$-free source; then send a point of $P$ to its fiberwise based map and descend. Conversely a section pulls back to a $G$-equivariant family of based maps and hence gives a map on $\Sigma^2P_+$.

\begin{proof}
    Since the $G$-action on $\Sigma^2P_+$ is free away from the base point, any $[f]\in \{\Sigma^2P_+,S^{k{\tilde{\R}}}\}^{G}$ can be represented by an $G$-map     \begin{equation}\label{equ:def-of-bf-from-smash-product}
    f: S^{\mathbb{R}^a} \wedge \Sigma^2P_+ \to S^{\mathbb{R}^a \oplus k{\tilde{\R}}}
        \end{equation}
    such that the $G$-action on ${\mathbb{R}^a}$ is trivial (see \cite[proposition 5.5]{Adams84}).
     Observe that 
     \[
     P\times_G(\mathbb R^{a+2})\cong\underline{\mathbb R}^{a+2} ,\qquad P\times_G(\mathbb R^a\oplus k{\tilde{\R}})
\cong\underline{\mathbb R}^{a}\oplus k\lambda
     \]
     are two bundles over $P\times_G * = Q$. The fiberwise compactification gives two (representation sphere-) bundles
     \[
     P\times_G S^{a+2}\cong\underline{S}^{a+2} ,\qquad P\times_G S^{(\mathbb R^a\oplus k{\tilde{\R}})}
\cong \prod_{q\in Q}S^{\mathbb R^a\oplus k\lambda_q} .
     \]
     
    Fix $p\in P$ and define
    \begin{align}
    f_p : S^{a+2} &\to S^{\mathbb{R}^a \oplus k{\tilde{\R}}}\label{equ:G-equiv-family-based-map} \\
    z &\mapsto f(z,p)\notag
    \end{align}
    which satisfies $f_{gp} = gf_p$. Hence it induces a well-defined bundle map 
    \begin{align}
P\times_G S^{a+2}\cong\underline{S}^{a+2} &\to P\times_G S^{(\mathbb R^a\oplus k{\tilde{\R}})}
\cong \prod_{q\in Q}S^{\mathbb R^a\oplus k\lambda_q}\label{equ:base-change-fiber}\\
[p,z] &\mapsto [p, f_p(z)]. \label{equ:bundle-map-of-Q}
    \end{align}
    over $P\times_G * = Q$. On the fiber over $q = [p] \in Q$ it is the map $f_p$ that does not depend on the choice of representative, so we may denote it by $f_{[p]}$. Note that the base point of $S^{\mathbb R^a}$ is $\infty$ and Lemma \ref{lem:based-map-adjunction} and (\ref{equ:def-of-bf-from-smash-product}) ensures $f_{[p]}(\infty) = f(\infty, p) = \infty$ for every $p\in P$. Thus $f_{[p]}$ is an element in desired $\operatorname{Map}_*(S^{a+2},S^{\mathbb R^a\oplus k\lambda_q})$. Notice that (\ref{equ:bundle-map-of-Q}) is continuous on $f$, so a homotopy in the $G$-equivariant mapping class corresponds to a homotopy of sections.

    Conversely, suppose we have a stable section
    \[
f\in \prod_{q\in Q}\mathcal E_k(q):= \prod_{q\in Q}\operatorname*{colim}_{a\to\infty}
\operatorname{Map}_*(S^{a+2},S^{\mathbb R^a\oplus k\lambda_q}).
    \]
    By the Freudenthal suspension theorem we can pick an universal $a\in \mathbb Z$ (that depends only on $k$) such that for every $q\in Q$,
    \[
    f_q \in \operatorname{Map}_*(S^{a+2},S^{\mathbb R^a\oplus k\lambda_q})
    \]
    is already in the stable range. 
    Then pull back this section along $P\to P\times_G * = Q$ we get a $G$-equivariant family of based maps as in (\ref{equ:G-equiv-family-based-map}). The rest of the proof explains the result of such pullback. By adjusting the first representative via the $G$-action, $f$ can be written as the form in (\ref{equ:bundle-map-of-Q}). Take $f_p$ there be the map induced by our $f$. We deduce that
    \begin{enumerate}
    \item According to the fiber of (\ref{equ:base-change-fiber}), we have $f_p\in\operatorname{Map}(S^{a+2},S^{\mathbb R^a\oplus k{\tilde{\R}}})$.
        \item $f$ is a based map on every fiber, so $f_p$ is a based map. We have $f_p\in\operatorname{Map}_*(S^{a+2},S^{\mathbb R^a\oplus k{\tilde{\R}}})$ for any $p\in P$. Again by Lemma \ref{lem:based-map-adjunction}, $f$ induces a map
        \[
         S^{\mathbb{R}^a} \wedge \Sigma^2(P_+) \to S^{\mathbb{R}^a \oplus k{\tilde{\R}}}
        \]
        sending $(z,p)$ to $f_p(z)$.
        \item From (\ref{equ:bundle-map-of-Q}) we have $[p,z] = [gp,gz] = [gp,z] \mapsto [gp, f_{gp}(z)] = [p, g^{-1}f_{gp}(z)]$, which indicates that $f_{gp} = gf_p$. Hence above map induced by $f$ is $G$-equivariant.
    \end{enumerate}
    Take sufficiently large $a$ (depends only on $P$ and $k$), we get a class in desired $\{\Sigma^2P_+,S^{k{\tilde{\R}}}\}^{G}$.

    Lastly, notice that the pullback is continuous on $f$, so a homotopy of sections corresponds to a homotopy in the $G$-equivariant mapping class.
\end{proof}

The fiber of $\mathcal E_k$ has homotopy group in degree $j\ge0$ equal to 
\begin{equation}\label{equ:homotopy-of-fiber}
\pi_j(\{S^2,S^{k\lambda_q}\}) \cong \pi^S_{j+2-k}.
\end{equation}
Thus for $k\ge3$ the fiber is $(k-3)$-connected. Its first nonzero group, in degree $k-2$, is $\pi^S_0=\mathbb Z$. The action of $\pi_1(Q)$ on this group is an isomorphism and the sign is given by the orientation character of $k\lambda$. We use $k=4n-2$, which is even, so the target bundle $\prod_{q\in Q}S^{\mathbb R^a\oplus k\lambda_q}$ is orientable. 
The canonical map from spheres to Eilenberg-MacLane sapces gives
\[
\operatorname{Th} (\underline{\mathbb R}^{a}\oplus k\lambda) \to K(\mathbb Z,a + k)
\]
which is fiberwise
\[
 S^{\mathbb R^a\oplus k\lambda_q}\longrightarrow K(\mathbb Z,a + k).
\]
It induces \[
\operatorname{Map}_*(S^{a+2},S^{\mathbb R^a\oplus k\lambda_q})\to \operatorname{Map}_*(S^{a+2},K(\mathbb Z,a + k)).
\]
\begin{definition}\label{def:primary-obstruction}
A section $\beta$ of $\mathcal E_k$ gives
\[
\beta : Q \to \prod_{q\in Q}\operatorname{Map}_*(S^{a+2},S^{\mathbb R^a\oplus k\lambda_q}) \to\operatorname{Map}_*(S^{a+2},K(\mathbb Z,a + k))\cong K(\mathbb Z, k-2).
\]
It defines a cohomology class 
\[
u(\beta)\in H^{k-2}(Q;\mathbb Z)
\]
which only depends on the homotopy class of sections. We call it the \textit{primary obstruction} of the stable section.
\end{definition}

Set $r=k-2=4n-4$. 
The first obstruction of the null-homotopy of the section $\beta: q \mapsto \mathcal E_k(q)$ is
\[
u(\beta)\in H^r(Q;\mathbb Z).
\]
The restriction of a section to an $r$-dimensional CW complex is classified by this primary class: the fibers are $(r-1)$-connected, so there are no lower obstructions; the choices on the $r$-cells are integers, and homotopies on the $(r-1)$-cells change them by cellular coboundaries. In particular, if $u(\beta)=0$, the restriction of $\beta$ to the $r$-skeleton is homotopic to the zero section. This argument concerns stable sections (sections of the stable mapping class bundle); on a finite CW complex a sufficiently large suspension $\Sigma^a$ (with trivial $G$-action) represents all the sections and homotopies in question.

\subsection{Obtain the family Seiberg-Witten from the stable section}
Family Seiberg-Witten is defined by counting the $S^1$-orbit. In family Bauer-Furuta setting, this corresponds to quotienting the domain by $S^1$-action, then pulling back the fundamental class of the target sphere. In this subsection, we show that corresponds to pulling back the primary obstruction of stable section to $S^1$-orbit, and consequently the characteristic hypothesis in Theorem \ref{thm:simply-connected-two-euler-annihilation} kills the primary obstruction.

\begin{prop}\label{prop:obtain-SW-from-primary-obstruction}
Define the natural projection
\[
p:S(n\mathbb H)/S^1\cong\mathbb{CP}^{2n-1}\longrightarrow Q = S(n\mathbb H)/G.
\]
 For an even integer $d$ and $\alpha\in\{S^{\mathbb R+n\mathbb H},S^{(4n-d){\tilde{\R}}}\}^{G}$, the corresponding stable section $\beta$ satisfies
\[
p^*u(\beta)=t_n(\operatorname{Res}_{S^1}^{G}\alpha).
\]
\end{prop}

\begin{proof}
Set $k=4n-d$. The quotient map
\[
p:S(n\mathbb H)/S^1\cong\mathbb{CP}^{2n-1}\longrightarrow Q = S(n\mathbb H)/G
\]
is a double cover. The pulled-back line bundle $p^*\lambda$ is canonically trivial: $S^1$-action is trivial on ${\tilde{\R}}$. Restricting $\alpha$ from $G$ to $S^1$ therefore corresponds to pulling back its section $\beta$ along $p$. The resulting ordinary stable map has integral cohomology class $p^*u(\beta)$ in degree $r=k-2$ (see the bottom-left square of the following commutative diagram).

The $S^1$-equivariant cofiber map identifies this map with the class defining $t_n(\operatorname{Res}\alpha)$. More explicitly, the groups $\{S^{\mathbb R},S^{k{\tilde{\R}}}\}^{S^1}$ and $\{S^{2\mathbb R},S^{k{\tilde{\R}}}\}^{S^1}$ vanish for $k\ge 3$, and quotienting the $S^1$-free source gives
\[
\begin{aligned}
\{S^{\mathbb R+n\mathbb H},S^{k{\tilde{\R}}}\}^{S^1}
&\cong \{\Sigma^2S(n\mathbb H)_+,S^{k{\tilde{\R}}}\}^{S^1}\\
&\cong \{\Sigma^2\mathbb{CP}^{2n-1}_+,S^k\}.
\end{aligned}
\]
These are the precise isomorphisms (15)--(16) used in \cite{Lin23}
. Then the pullback of $1\in H^k(S^k)$ defines $t_n$. This definition and the suspension property of $H^*$ give the bottom-right square of the following commutative diagram:

\begin{tikzcd}[
  column sep=2em,
  row sep=3.2em,
  cells={nodes={inner sep=2pt}}
]
&[-1.5em] \alpha \arrow[r, mapsto]
& \Res^G_{S^1}\alpha
&
\\[-2.7em]
& \left\{S^{\R+n\HH},S^{k{\tilde{\R}}}\right\}^{G}
    \arrow[r, "\Res"]
    \arrow[d, "\cong"',"c_n", description]
& \left\{S^{\R+n\HH},S^{k{\tilde{\R}}}\right\}^{S^1}
    \arrow[d, "\cong"]
&
\\
& \left\{\Sigma^2S(n\HH)_+,S^{k{\tilde{\R}}}\right\}^{G}
    \arrow[d, "\cong"']
& \left\{\Sigma^2S(n\HH)_+,S^{k{\tilde{\R}}}\right\}^{S^1}
    \arrow[d, "\cong"]
&
\\
\beta\in
    \arrow[d, mapsto]
& \left\{
      Q,\prod\limits_{q:Q}\{
        S^{2},S^{ k\lambda_q}
      \}
    \right\}
    \arrow[d,"u", description]
    \arrow[dr, "p^*"']
& \left\{\Sigma^2\mathbb{CP}^{2n-1}_+,S^{k}\right\}
    \arrow[r, "t_n"]
    \arrow[d, "\cong", description]
    \arrow[dr, phantom, "\circlearrowright", description]
& H^k\!\left(\Sigma^2\mathbb{CP}^{2n-1}_+;\Z\right)
    \arrow[d, "\cong", description]
\\
u(\beta)\in
& \left\{Q;K(\Z,k-2)\right\}
    \arrow[dr, "p^*"']
& \{\mathbb{CP}^{2n-1},\{S^2,S^k\}\}
    \arrow[d]
& H^{k-2}\!\left(\mathbb{CP}^{2n-1};\Z\right)
\\
&
& \left\{\mathbb{CP}^{2n-1},K(\Z,k-2)\right\}
    \arrow[ur, "\cong"']
&
\end{tikzcd}

Chasing this diagram gives  
\[t_n(\operatorname{Res}_{S^1}^{G}\alpha) = p^*u(\beta).
\]
\end{proof}

Therefore the hypothesis $t_n(\operatorname{Res}_{S^1}^{G}\alpha)=0$ says $p^*u(\beta)=0$.

\begin{lemma}\label{lem:SW-zero-means-primary-obstruction-zero}
    If $\alpha\in\{S^{\mathbb R+n\mathbb H},S^{(4n-2){\tilde{\R}}}\}^{G}$ satisfies
\[
t_n(\operatorname{Res}_{S^1}^{G}\alpha)=0,
\]
then the corresponding stable section $\beta$ satisfies
\[
u(\beta)=0.
\]
\end{lemma}
\begin{proof}
We check injectivity of $p^*$ in this degree. The diagonal left action of unit quaternion $\mathrm{Sp}(1)$ on $S(n\mathbb H)$ gives a fibration
\[
\mathbb{RP}^2=\mathrm{Sp}(1)/G\longrightarrow Q
\longrightarrow \mathbb{HP}^{n-1}.
\]
The base is simply connected. Its integral cohomology is $\mathbb Z$ in degrees $0,4,\ldots,4n-4$ and zero otherwise. The fiber has integral cohomology $\mathbb Z$ in degree zero and $\mathbb Z/2$ in degree two, and no other positive cohomology. In the integral cohomology Serre spectral sequence there are only rows zero and two; the only possible nonzero differential between them is $d_3$, whose target base degree is $3$ modulo $4$ and hence zero:
\begin{center}
\begin{tikzpicture}[
    x=1.6cm,
    y=1.1cm,
    every node/.style={font=\large}
]

\draw[thick] (0,0) -- (7.3,0);
\draw[thick] (0,0) -- (0,3.5);

\node[above left] at (0,3.5) {$\mathbb{RP}^{2}$};

\node[left] at (0,0.5) {$0$};
\node[left] at (0,1.5) {$1$};
\node[left] at (0,2.5) {$2$};

\foreach \y in {0.5,1.5,2.5} {
    \node at (0.7,\y) {$\cdots$};
    \node at (2.0,\y) {$\mathbb{Z}_2$};
    \node at (3.1,\y) {$0$};
    \node at (4.2,\y) {$0$};
    \node at (5.3,\y) {$0$};
    \node at (6.4,\y) {$\mathbb{Z}_2$};
}

\node[fill=white,inner sep=3pt] at (2.0,0.5) {$\mathbb{Z}$};
\node[fill=white,inner sep=3pt] at (2.0,1.5) {$0$};
\node[fill=white,inner sep=3pt] at (6.4,0.5) {$\mathbb{Z}$};
\node[fill=white,inner sep=3pt] at (6.4,1.5) {$0$};

\node[below] at (2.0,0) {$4n-8$};
\node[below] at (6.4,0) {$4n-4$};

\node[right] at (7.3,0) {$\mathbb{HP}^{n-1}$};

\end{tikzpicture}
\end{center}
Therefore the spectral sequence collapses. In total degree $r=4n-4$ its only nonzero term is the $\mathbb Z$ in bidegree $(r,0)$. It follows that $H^r(Q;\mathbb Z)=\mathbb Z$, with no torsion. The integral transfer for a double cover 
\[
H^r(\mathbb{CP}^{2n-1};\mathbb Z)
\xrightarrow{\,\operatorname{tr}\,}
H^r(Q;\mathbb Z).
\]
is defined by
\[
(\operatorname{tr}\varphi)(\tau)
=\varphi(\widetilde\tau_1)+\varphi(\widetilde\tau_2).
\]
for any simplex $\tau$ in $Q$ and any cochian $\varphi$ of the covering space $\mathbb{CP}^{2n-1}$, where $\widetilde\tau_i$ are lifts of $\tau$. On the other hand, the pullback
\begin{equation}\label{equ:pullback-homology-at-r}
H^r(Q;\mathbb Z)
\xrightarrow{\,p^*\,}
H^r(\mathbb{CP}^{2n-1};\mathbb Z)
\end{equation}
for any cochain $c$ of $Q$, satisfies
\[
p^*c(\widetilde\tau_i) = c(\tau).
\]
Hence
\[
(\operatorname{tr}\,p^*c)(\tau) = p^*c(\widetilde\tau_1)+p^*c(\widetilde\tau_2) = 2c(\tau)
\]
Hence the integral transfer for a double cover satisfies
\[
\operatorname{tr}\,p^*=2\operatorname{id}.
\]
Hence $p^*$ is injective in degree $r=4n-4$. We conclude that $u(\beta)=0$.
\end{proof}

\section{Applications}
\subsection{Lin-Mukherjee vanishing theorem}
Lin-Mukherjee vanishing theorem now can be deduced from the stable section construction easily:
\begin{proof}[Proof of Theorem \ref{thm:LM-vanishing}]
    When $2\chi(X) + 3\sigma(X) > 4$, the setting
    \[
    n=-\frac{\sigma}{16}
    \]
    gives that
    \[
    b^+ > 4n.
    \]
    Hence the connectivity of the fiber (\ref{equ:homotopy-of-fiber}) is greater than $4n-2$. However, the $\operatorname{Pin}(2)$-orbit has 
    \[
    \operatorname{dim}(S(n\mathbb H)/\operatorname{Pin}(2))=4n-2.
    \]
Hence any section is homotopic to the trivial section, even though the bundle $k\lambda$ might be non-orientable. Therefore, the invariants
$\operatorname{FBF}^{\operatorname{Pin}(2)}(X, f, \widetilde{f} )$ vanishes.

When $2\chi(X) + 3\sigma(X) = 4$, we don't have such strong argument: the stable section does have the primary obstruction 
\[
u(\beta) \in H^{4n-2}(S(n\mathbb H)/\operatorname{Pin}(2);\mathbb Z)\cong \mathbb Z_2.
\]
However, we can apply Proposition \ref{prop:G-equiv-to-stable-section} to \[
G= S^1 ,\qquad Q = S(n\mathbb H)/S^1 = \mathbb{CP}^{2n-1}.
\]
The primary obstruction now is 
\[
u(\beta) \in H^{4n-2}(S(n\mathbb H)/S^1;\mathbb Z)\cong \mathbb Z.
\]
From the commutative diagram in Proposition \ref{prop:obtain-SW-from-primary-obstruction}, we know this primary obstruction is a pullback from $H^{4n-2}(S(n\mathbb H)/\operatorname{Pin}(2);\mathbb Z)\cong \mathbb Z_2$. This pullback homomorphism $\mathbb Z_2 \to \mathbb Z$ must be trivial. As $\operatorname{dim}(\mathbb{CP}^{2n-1}) = 4n-2$, this primary obstruction is the only obstruction. Therefore, the invariants
\[
\operatorname{{Res}_{S^1}}\operatorname{FBF}^{\operatorname{Pin}(2)}(X, f, \widetilde{f} ) = \operatorname{FBF}^{S^1}(X, f, \widetilde{f} )\]
vanishes.
\end{proof}

\subsection{Vanishing theorem after two stabilizations}
We first explain why two Euler classes kill the $G$-equivariant map.

\begin{proof}[Proof of Theorem \ref{thm:simply-connected-two-euler-annihilation}]
Set
\[
P=S(n\mathbb H),\qquad Q=S(n\mathbb H)/G,\qquad\lambda=S(n\mathbb H)\times_G{\tilde{\R}},
\]
where $\lambda$ is the associated real line bundle on $Q$. The action on $S(n\mathbb H)$ is free: if a unit quaternion fixes a nonzero component of a vector, that quaternion is $1$. The cofiber sequence
\[
S^{\mathbb R}\xrightarrow{i_n} S^{\mathbb R+n\mathbb H}\xrightarrow{c_n}
\Sigma^2S(n\mathbb H)_+
\]
and exactness give a class
\[
\beta\in\{\Sigma^2S(n\mathbb H)_+,S^{(4n-2){\tilde{\R}}}\}^{G}
\]
such that
\[
\alpha =\beta c_n.
\]
The Euler class $e \in \{S^0, S^{\tilde{\R}}\}^G$ corresponds to the zero section of 
\[
\prod_{q\in Q}\{S^0, S^{\lambda_q}\}.
\]
We call this section the stable Euler section $E(\lambda)$. The product of the $G$-equivariant class corresponds to the direct sum of representations (smash product of representation spheres). Hence $e^2$ corresponds to the stable Euler section $E(2\lambda)$. Note that there are two different sections:
\begin{enumerate}
    \item The zero section of the vector bundle $q\mapsto 0_q$, which represents the Euler class;
    \item The base-point section $q\mapsto \infty_q$, which is the trivial section and corresponds to trivial $G$-equivariant map.
\end{enumerate}

The smooth compact manifold $Q$ has dimension $4n-2=r+2$
. By Lemma \ref{lem:SW-zero-means-primary-obstruction-zero} and the obstruction calculation around Definition \ref{def:primary-obstruction}, $\beta$ may be represented by a stable section which is trivial on $Q^{(r)}$, using homotopy extension property (HEP) of $Q$. The class $e^2$, over $Q$, is the stable Euler section $E(2\lambda)$ of the one-point compactification of $2\lambda$. Its restriction to $Q^{(1)}$ is homotopy to the trivial section: every rank-two vector bundle over a one-dimensional CW complex has a nowhere zero section, say $v(q)$, obtained by choosing unit vectors on vertices and extending along edges because the unit circle is connected. Then 
\[
H(q,t)=
\begin{cases}
\dfrac{t}{1-t}\,v(q),&0\le t<1,\\[4pt]
\infty_q,&t=1.
\end{cases}
\]
supplies a homotopy from the zero Euler section to the section at infinity. Thus by HEP we can choose a representative for $E(2\lambda)$ that equals the base-point section on $Q^{(1)}$.

The multiplication of these classes corresponds to direct sum of representations, which can be obtained by exterior product of representations (fiberwise exterior smash of representation spheres) over $Q\times Q$:
\[
\beta e^2=\Delta^*(\beta\boxtimes e^2),
\qquad
\Delta(q)=(q,q),
\]
where
\[
(\beta\boxtimes e^2)(q_1,q_2)
=\beta(q_1)\wedge E(2\lambda)(q_2).
\]
As the smash product is the base point if one of these sections is at the base point, $\beta\boxtimes e^2$ is the base point on
\[
(Q^{(r)}\times Q)\ \cup\ (Q\times Q^{(1)}).
\]
Every product cell of total dimension at most $r+2$ is in that union: a cell outside the union would have first dimension at least $r+1$ and second dimension at least $2$, with sum at least $r+3$. Hence 
\[
(Q\times Q)^{(r+2)} \subset (Q^{(r)}\times Q)\ \cup\ (Q\times Q^{(1)}).
\]
But $Q$ is $(r+2)$-dimensional. By the cellular approximation theorem we have a homotopy from $\Delta$ to an approximated diagonal
\[
\Delta'(Q)\subset (Q\times Q)^{(r+2)}.
\]
The approximated diagonal therefore pulls the exterior section back to the trivial section. Pullbacks along homotopic diagonals give homotopic stable sections after the induced bundle identification. Along the actual diagonal $\Delta$ the target virtual bundle is
\[
((4n-2)\lambda-\mathbb R^2)\oplus2\lambda
=4n\lambda-\mathbb R^2.
\]
The resulting section is exactly the class $\beta e^2$. Consequently
\[
\beta e^2=0,
\qquad
\alpha e^2=(\beta e^2)c_n=0.
\]
This completes the assertion without a choice-dependent claim about the lift.
\end{proof}

Now we prove Corollary \ref{coro:main-theorem-2stabilization-vanish}. Before that we show that under our convention, the parameter circle is mapped to the base point.
\begin{lemma}\label{lem:alpha-i-trival}
     Let X be an oriented spin 4-manifold with $b_1(X) = 0$. We fix a homology orientation and
a spin structure on X. Let $f : X \to X$ be a diffeomorphism which preserves the spin structure and acts
trivially on $H^*(X; \mathbb R)$. Suppose $2\chi(X) + 3\sigma(X) \ge -4$. Let $ n=-\frac{\sigma}{16}$ and let $i_n:S^{\mathbb R}\to S^{\mathbb R+n\mathbb H}$ be inclusion at zero in the quaternionic coordinates. For any lift $\widetilde{f}$ of $f$ to an automorphism on the spin structure, 
\[\alpha:=\operatorname{FBF}^{\operatorname{Pin}(2)}(X, f\# \operatorname{id}_{2(S^2\times S^2)}, \widetilde{f} )\]
satisfies
    \[
    \alpha\circ i_n=0.\]
\end{lemma}
\begin{proof}
    By the assumption $2\chi(X) + 3\sigma(X)$ can only be $-4$ or nonnegative, but the latter case would make this invariant of $X\#2(S^2\times S^2)$ be at least $8$ and solved by case (2) of Theorem \ref{thm:LM-vanishing}. For simplicity we only consider $2\chi(X) + 3\sigma(X) = -4$.

    From $\chi(X) = 2+b^++b^-$ we solve that 
    \begin{align*}
        b^+ &= 4n-2,\\
        b^- &= 20n-2.
    \end{align*}
    As $b^+ \ge 0$ and $n$ is integer by Rokhlin Theorem, $n\ge 1$.

    Choose a finite-dimensional representative of $\alpha\circ i_n$. For sufficiently large $a,c,k$, the finite-dimensional approximation of the family BF map gives a representative
\[
\alpha:
S^{(a+1)\mathbb R+c{\tilde{\R}}+(k+n)\mathbb H}
\longrightarrow
S^{a\mathbb R+(c+b^+){\tilde{\R}}+k\mathbb H}.
\]
Setting the last $n$ quaternionic coordinates equal to zero gives a representative of $x\circ i_n$,
\begin{equation}\label{equ:bf-map-without-spinor}
F:
S^{(a+1)\mathbb R+c{\tilde{\R}}+k\mathbb H}
\longrightarrow
S^{a\mathbb R+(c+b^+){\tilde{\R}}+k\mathbb H}.
\end{equation}
The source still contains $k\mathbb H$. We first construct a nullhomotopy where these remaining quaternionic coordinates vanish.

Now construct the nullhomotopy on the $S^1$-fixed subspace. The $S^1$-fixed subspace of the source of (\ref{equ:bf-map-without-spinor}) is $S^{(a+1)\mathbb R+c{\tilde{\R}}}$. On the mapping torus, taking gauge-$S^1$ fixed points sets the spinor equal to zero. On each fiber $X_t$, the monopole map then reduces to
\[
d_t^*\oplus d_t^+:
\Omega^1(X_t;i\mathbb R)
\longrightarrow
\bigl(\Omega^0(X_t;i\mathbb R)/i\mathbb R\bigr)
\oplus\Omega^+(X_t;i\mathbb R).
\]
Its kernel is zero because $b_1(X)=0$, and its cokernel is the space of self-dual harmonic forms on $X_t$:
\[
\begin{array}{ccc}
U_t^+ & \xrightarrow{\ d_t^*\oplus d_t^+\ } &
(d_t^*\oplus d_t^+)(U_t^+)\oplus\mathcal H^+(X_t;i\mathbb R)\\[3pt]
\big\downarrow{\cong} && \big\downarrow{\cong}\\[3pt]
\mathbb R^c & \xrightarrow{\ u\mapsto(u,0)\ } &
\mathbb R^c\oplus\mathbb R^b 
\end{array}
\]
Here $U^+_t \subset \Omega^1(X_t;i\mathbb R)$ is an $c$-dimensional approximation. Since $f^*$ acts trivially on $H^2(X;\mathbb R)$, the compatible trivializations in the family Seiberg-Witten construction identify the finite-dimensional approximation of this operator with the same coordinate inclusion $u\mapsto(u,0)$ for every $t$.

Consequently, in the reduced family BF class, the parameter circle in the $S^1$-fixed subspace factors through
\[
S^{2\mathbb R}
\xrightarrow{\,c_0\,}
S^{\mathbb R}\wedge S(2\mathbb R)_+
\xrightarrow{\,1\wedge p\,}
S^{\mathbb R}.
\]
Here $c_0$ is the Pontryagin–Thom collapse for the parameter circle, and $p:S(2\mathbb R)_+\to S^0$ sends the whole parameter circle to the non-base point. Topologically $c_0$ identify the north and the south pole; The result looks like a pumpkin, with top and bottom glued together. Then each slice of the pumpkin is mapped to the circle identically:
\begin{center} 
\begin{tikzpicture}[
    x=1.15cm, y=1.15cm,
    line width=.85pt, line cap=round, line join=round,
    hidden/.style={dash pattern=on 3pt off 3pt},
    map/.style={-{Stealth[length=5pt,width=5pt]}}
]

\draw (0,0) circle (1);
\draw[hidden] (-.99,-.14)
    arc[start angle=180,end angle=0,x radius=.99,y radius=.34];
\draw (-.99,-.14)
    arc[start angle=180,end angle=360,x radius=.99,y radius=.40];

\draw[map] (1.65,0) -- node[above=3pt] {$c_0$} (2.95,0);

\begin{scope}[shift={(4.82,0)}]
  \coordinate (P) at (0,.66);       
  \coordinate (B) at (-.09,-.56);   
  \coordinate (L) at (-.99,.64);
  \coordinate (R) at (1.00,.65);
  \coordinate (U) at (-.54,1.24);
  \coordinate (V) at (.54,1.30);

  \draw (L) .. controls (-1.00,.91) and (-.76,1.18) .. (U)
            .. controls (-.22,1.42) and (.43,1.46) .. (V)
            .. controls (.97,1.36) and (1.08,1.03) .. (R)
            .. controls (1.48,.22) and (.91,-.82) .. (B)
            .. controls (-.88,-.60) and (-1.25,-.34) .. (-1.09,.39)
            .. controls (-1.08,.53) and (-1.03,.59) .. cycle;

  \draw[hidden] (P)
      .. controls (-.09,.28) and (-.31,-.44)
      .. (-.59,-.53);
  \draw[hidden] (.06,.39)
      .. controls (.26,-.49) and (.51,-.78) .. (.99,.43);

  \draw (U) .. controls (-.11,1.30) and (.15,1.04) .. (P)
            .. controls (.11,.94) and (.28,1.34) .. (V);
  \draw (L) .. controls (-.65,1.00) and (-.12,.97) .. (P)
            .. controls (.33,.80) and (.66,.84) .. (R);
  \draw (P) .. controls (-.35,.62) and (-.85,.10) ..
  (-.59,-.53);
  \draw (P) .. controls (.59,.38) and (.68,-.20) .. (B);

  \fill (P) circle[radius=1.35pt];
\end{scope}

\draw[map] (6.65,0) -- node[above=3pt] {$1\wedge p$} (7.74,0);

\draw (9.02,0) circle[radius=.54];
\fill (9.02,0) ++(32:.54) circle[radius=1.65pt];
\end{tikzpicture}
\end{center}
This composite is nullhomotopic because $\pi_2(S^1)=0$. The residual $C_2$-action on the parameter spheres is trivial, so the nullhomotopy is $C_2$-equivariant. After increasing $a,c$ if necessary, we obtain an actual equivariant nullhomotopy of $F$ on $S^{(a+1)\mathbb R+c{\tilde{\R}}}$. The inclusion of this subspace is a $G$-cofibration. By the equivariant version of homotopy extension property (HEP), we may therefore assume that $F$ sends the entire subspace to $\infty$.

Now we can perturb to remove the zeros on the $G$-free part. 
Outside $S^{(a+1)\mathbb R+c{\tilde{\R}}}$, the $G$-action there is free: if $gq=q$ for a nonzero quaternion $q$, then $g=1$. Hence the free quotient has dimension
\[
\dim\!\left(
\bigl(S^{(a+1)\mathbb R+c{\tilde{\R}}+k\mathbb H}
\setminus S^{(a+1)\mathbb R+c{\tilde{\R}}}\bigr)/G
\right)
=(a+1+c+4k)-1=a+c+4k.
\]
Over this quotient, $F$ defines a section of the associated sphere bundle with fiber $S^{a\mathbb R+(c+b^+){\tilde{\R}}+k\mathbb H}$. The origin in each fiber defines a zero section of codimension $
a+c+b^++4k$. Since $F$ equals $\infty$ on the $S^1$-fixed subspace, all its zeros lie in a compact subset of the $G$-free part. We may thus apply ordinary relative transversality on the quotient, with the perturbation supported away from the $S^1$-fixed subspace. The transverse zero set has expected dimension
\[
(a+c+4k)-(a+c+4k+b^+)=-b^+<0,
\]
so it is empty.

Lifting this perturbation gives a $G$-equivariant map homotopic to $F$ whose image avoids the origin. The complement of the origin in a representation sphere contracts $G$-equivariantly to $\infty$. Thus $F$ is $G$-equivariantly nullhomotopic, proving $\alpha\circ i_n=0$. This finite-dimensional nullhomotopy also gives the asserted vanishing in a complete universe.
\end{proof}

\begin{proof}[Proof of Corollary \ref{coro:main-theorem-2stabilization-vanish}]
  Again we only deal with $2\chi(X) + 3\sigma(X) = -4$. 
  Set 
  \[
  \alpha = \operatorname{FBF}^{\operatorname{Pin}(2)}(X, f\# \operatorname{id}_{2(S^2\times S^2)}, \widetilde{f} ).
  \]
  Above lemma shows 
  \[
  \alpha\circ i_n = 0.
  \]
  The family Seiberg-Witten invariant is given by $t_n\circ\operatorname{Res}_{S^1}$ (see \cite{Lin23} Section 2.1.2). So the hypothesis on family Seiberg-Witten invariant shows $t_n(\operatorname{Res}_{S^1}(\alpha))= 0 $. By Proposition \ref{prop:G-equiv-to-stable-section}, we have
  \[
  \alpha \cdot e^2 = 0.
  \]
  Corollary 2.27 of \cite{Lin23} states that, fiberwise connected sum with either product or twisted spin family on $(S^2\times S^2)\times S^1$, multiplies the $\operatorname{Pin}(2)$-FBF by $e$. Iterating it identifies $\alpha e^2$ with the invariants after the two identity stabilizations.
\end{proof}

\bibliographystyle{alpha}
\bibliography{./vanishing}
\end{document}